\documentclass[12pt]{article}
\renewcommand{\baselinestretch}{1.150}

\usepackage{amsmath}
\usepackage{amsthm}
\usepackage{amssymb, amsfonts}
\usepackage{mathrsfs}
\newcommand{\ds}{\displaystyle}
\newcommand{\p}{\partial}

\newtheorem{prop}{Proposition}

\newtheorem{conj}{Conjecture}

\DeclareMathAlphabet{\bi}{OML}{cmm}{b}{it}

\newcommand{\be}{\begin{equation}}
\newcommand{\ee}{\end{equation}}
\newcommand{\ba}{\begin{array}}
\newcommand{\ea}{\end{array}}
\def\dder#1{\p_{#1}}

\date{March 1, 2017}

\title{\protect\vspace{-2cm}\sc\Large A Simple Construction of Recursion Operators\\ for Multidimensional Dispersionless Integrable Systems}
\author{{\sc A. Sergyeyev}\\[5mm]
Mathematical Institute, Silesian University in Opava,\\ Na
Rybn\'\i{}\v{c}ku 1, 746 01 Opava,~Czech Republic\\
E-mail: \tt{Artur.Sergyeyev@math.slu.cz}}

\begin{document}
\maketitle

\begin{abstract}\protect\vspace*{-10mm}
We present a simple novel construction of recursion operators for
integrable multidimensional dispersionless systems that admit a
Lax pair whose operators are linear in the spectral
parameter and do not involve the derivatives with respect to the
latter. New examples of recursion operators
obtained using our technique include {\em inter alia} those for the
general heavenly equation, which describes a class of anti-self-dual solutions of the vacuum Einstein equations, and a six-dimensional equation resulting
from a system of Ferapontov and Khusnutdinova.
\end{abstract}


\section*{Introduction}
The field of integrable systems is a dynamic subject
closely related to many parts of
modern mathematics and featuring manifold applications, cf.\ e.g.\ \cite{ac, a-s, bl, blu, ca, cgs, dun, l, n, nd, olverbook, vps}
and references therein.\looseness=-1

A characteristic feature of integrable systems is that they never show up alone -- instead, they
always belong to integrable hierarchies, and the latter are typically
constructed using the recursion operators (ROs), cf.\ e.g.\ \cite{ac, bl, dun, ma, n, olver, olverbook} and references therein.
In particular, the ROs for multidimensional integrable
dispersionless systems
are now a subject of intense research, cf.\ e.g.\
\cite{ac, dun, m-s, ms3, ms-ip, m1,m2,m3, olverbook} and references therein. In fact, the overwhelming majority of integrable systems in four or more independent variable known to date, including those relevant for applications, like the (anti-)self-dual Yang--Mills equations or the (anti-)self-dual vacuum Einstein equations, is dispersionless, i.e., can be rewritten as first-order homogeneous quasilinear systems, cf.\ e.g.\ \cite{ac, dun, s} and Example~5 below for details.\looseness=-1

The study of a hierarchy associated to a given integrable system, and hence of its recursion operator,
is important for many reasons, {\em inter alia} because known exact solutions of integrable systems, like, say, multisoliton and finite-gap solutions, always happen to be invariant under a certain combination of the flows in the hierarchy, and this is instrumental in the actual construction of the solutions in question, cf.\ e.g.\ \cite{dun, n}.\looseness=-1

Another feature making recursion operators relevant for applications is the fact that existence of a (formal) recursion operator can, 
under certain additional assumptions, be employed as an efficient integrability test or as a means of proving nonintegrability, see e.g.\ \cite{my} and references therein for general theory and \cite{sv} for a recent example.\looseness=-1


There are several methods for construction of the ROs for the
integrable dispersionless systems: the partner symmetry approach \cite{m-s}, the
method based on adjoint action of the Lax operators \cite{ms-ip},
and the approach using the Cartan equivalence method \cite{m1, m2,
m3}.

Below we present another method for construction of ROs for
integrable multidimensional dispersionless systems. Roughly
speaking, it consists in constructing a Lax pair for
the system under study with the following additional properties: i)
the Lax operators are linear in the spectral parameter and ii) the
solution of the associated linear system is a (nonlocal) symmetry of
the nonlinear system under study. Taking the coefficients at the
powers of spectral parameter in the operators constituting the newly
constructed Lax pair yields, under certain technical
conditions, the RO. A similar approach, with symmetries replaced by
cosymmetries and Proposition~\ref{Th1} by Proposition~\ref{Th2},
allows one to look for the adjoint ROs.\looseness=-1

To the best of author's knowledge, the key idea of our approach
that the Lax pair employed for the construction
of the recursion operator
should be {\em constructed} from the original Lax pair but is by no means obliged to be identical with the
latter, has not yet appeared in the literature. This kind of flexibility
noticeably enhances the applicability of our method.
On the other hand,
the idea of extraction of a recursion operator from a Lax pair which is
linear in the spectral parameter was already explored 
e.g.\ in \cite{mns, p}.\looseness=-1


The paper is organized as follows. In Section~\ref{pr} we recall
some basic facts\nopagebreak[1] concerning the geometric theory of PDEs. In
Section~\ref{rk} we state our results on the construction
of recursion operators and adjoint recursion operators. Section~3
explores the construction of special Lax pairs that
give rise to ROs through the results of Section~\ref{rk}.
Section~\ref{ex} gives a selection of examples illustrating our
approach, and Section~\ref{co} provides a brief discussion.\looseness=-1

\section{Preliminaries}\label{pr}
Let $m,N,d$ be nonzero natural numbers. Consider a system of $m$ PDEs
\begin{equation}\label{sys}F_I=0,\quad I=1,\dots,m,\end{equation}
in
$d$ independent variables $x^i$, $i=1,\dots,d$, for an un\-known $N$-com\-ponent vector
function $\bi{u}=(u^1,\dots,u^N)^{T}$.
Here and below the superscript `$T$' indicates the transposed matrix. For the sake of brevity we shall
employ the notation $\vec
x=(x^1,\dots,x^d)$ and $\mathcal{F}=(F_1,\dots,F_m)^T$.

As usual, we put $u^\alpha_{i_1\dots i_n}=\p^{i_1+\cdots+i_n}u^\alpha/\p
(x^1)^{i_1}\cdots \p  (x^n)^{i_n}$ and $u^\alpha_{00\cdots0}\equiv
u^\alpha$.

Following e.g.\ \cite{olverbook, kv}, $x^i$ and $u^\alpha_{i_1\dots
i_n}$ are considered here and below as independent quantities and can be
viewed as coordinates on an abstract infinite-dimensional space (a
jet space).
By a {\it local function} 
we shall mean a 
function of $\vec x$, $\bi u$ and of a finite number of
the derivatives of $\bi{u}$. To indicate locality of a function $f$
we shall write it as $f(\vec x,[\boldsymbol{u}])$. In what follows all functions
are tacitly assumed to be locally meromorphic
in 
all their arguments.

We denote the total derivatives with respect to $x^j$ by
\[
{\displaystyle D_{x^j} = \frac{\p}{\p x^j}
 } + \displaystyle\sum\limits_{\alpha=1}^{N}
 \sum\limits_{i_1,\dots,i_n=0}^{\infty}{\displaystyle
   u^\alpha_{i_1\dots i_{j-1} (i_j+1) i_{j+1}\dots i_n} \frac{\p}{\p u^\alpha_{i_1\dots i_n}}},
\]
cf.\ e.g.\ \cite{olverbook, kv}. For the sake of simplicity
the extensions of the above $D_{x^i}$ to nonlocal variables are again denoted by $D_{x^i}$, as e.g.\ in
\cite{ms-ip}.\looseness=-1

The system (\ref{sys}) together with all its differential consequences
\[D_{x^1}^{i_1}\cdots D_{x^n}^{i_n} F_I = 0,\quad i_1,i_2,\dots,i_n=1,\dots,d,\quad n=1,2,\dots,\quad I=1,\dots,m,
\]
defines a {\em diffiety} $\mathrm{Sol}_{\mathcal{F}}$, a submanifold of the infinite jet space,
see e.g.\ \cite{kv} and references therein for details.
This submanifold can be informally thought of as a  set of all formal solutions of (\ref{sys}),
which motivates the choice of notation. 
Below all equations will be required to hold only on
$\mathrm{Sol}_\mathcal{F}$, or a certain differential covering
thereof, see e.g.\ \cite{kv} for details on the latter, rather than e.g.\ on the whole jet space, unless
otherwise explicitly stated.
\looseness=-1

The {\it directional derivative}, cf.\ e.g.\ \cite{bl, ms-ip} and references therein, along 
$\bi{U}=(U^1,\dots,U^N)^{T}$
is the following vector field on the jet space
$$
\dder{\bi{U}} = \sum_{\alpha=1}^{N}\sum_{i_1,\dots,i_n=0}^{\infty}
(D_{x^1}^{i_1}\cdots D_{x^n}^{i_n} U^\alpha)
\frac{\p}{\p u^\alpha_{i_1\dots i_n}}.
$$
The total derivatives as well as
the directional derivative can be applied to (possibly vector or matrix) local functions
$P$.

Recall, cf.\ e.g.\ \cite{olverbook, kv}, that
a local $N$-component vector function $\bi{U}$ is a (characteristic
of a) {\it symmetry} for the system (\ref{sys}) if $\bi{U}$
satisfies the linearized version of this system, namely,
$\ell_{\mathcal{F}}(\bi U)= 0$, where
\[
\ell_f=\sum_{\alpha=1}^{N}\sum_{i_1,\dots,i_n=0}^{\infty}
\frac{\p f}{\p u^\alpha_{i_1\dots i_n}} D_{x^1}^{i_1}\cdots D_{x^n}^{i_n}
\]
is the operator of linearization and
$\ell_{\mathcal{F}}=(\ell_{F_1},\dots,\ell_{F_m})^T$.

We stress that by our
blanket assumption the condition $\ell_{\mathcal{F}}({\bi U})=0$ is required to hold on
$\mathrm{Sol}_\mathcal{F}$ only rather than on the whole jet space.
Informally 
this condition means that the flow
$\bi{u}_\tau=\bi{U}$ leaves $\mathrm{Sol}_\mathcal{F}$ invariant.\looseness=-1

Notice an
important identity
$\ell_{f}(\boldsymbol{U})=\partial_{\boldsymbol{U}}(f)$
which holds for any local $f$ and $\boldsymbol{U}$, see e.g.\ \cite{kv}.

Finally, recall that
if we have an operator in total derivatives
\[
Q=Q^0+\sum\limits_{k=1}^q \sum\limits_{j_1=1}^d\cdots
\sum\limits_{j_k=1}^d Q^{j_1\dots j_k} D_{x^{j_1}}\cdots
D_{x^{j_k}},
\]
where $Q^0$ and $Q^{j_1\dots j_k}$ are $s\times s'$-matrix-valued
local functions, its formal adjoint reads (cf.\ e.g.\ \cite{olverbook})
\[
Q^\dagger=(Q^0)^T+\sum\limits_{k=1}^q \sum\limits_{j_1=1}^d\cdots
\sum\limits_{j_k=1}^d (-1)^k D_{x^{j_1}}\cdots D_{x^{j_k}}\circ
(Q^{j_1\dots j_k})^T.
\]


For integrable systems in more than two independent variables their
symmetries usually depend, in addition to local variables, described
above, also on nonlocal variables. In other words, we are forced to
consider solutions $\boldsymbol{U}$ of
$\ell_{\mathcal{F}}(\boldsymbol{U})=0$ involving nonlocal variables.
Some authors refers to such objects as to the shadows of nonlocal
symmetries, see e.g.\ \cite{kv} and references therein for more
details on this terminology; cf.\ e.g.\ also \cite{blu}.
However, as in the present paper we shall not deal with full nonlocal symmetries
in the sense of \cite{kv}, for the sake of simplicity in
what follows we shall refer to solutions $\boldsymbol{U}$ of
$\ell_{\mathcal{F}}(\boldsymbol{U})=0$ involving nonlocal variables
just as to {\em nonlocal symmetries}.\looseness=-1

A {\em cosymmetry} $\boldsymbol{\gamma}$, also known as an adjoint symmetry, see e.g.\ \cite{blu},
is a quantity which is dual to a symmetry: it has $m$ rather than $N$ components and
satisfies (cf.\ e.g.\ \cite{bl,kv} for details) 
the 
system
$\ell_{\mathcal{F}}^{\dagger}(\boldsymbol{\gamma})=0$.
Note that cosymmetries, just like symmetries, may depend on nonlocal variables in addition to the local ones.

A {\em recursion operator} (RO) for our system (\ref{sys}) is then, roughly speaking, an operator that sends a symmetry of (\ref{sys}) into
another symmetry thereof \cite{olver, olverbook}. However, this intuitive definition works well only when the RO contains no nonlocal terms, as often is the case for linear and linearizable equations, cf.\ e.g.\ \cite{bp, gut, marv, olverbook} and references therein.\looseness=-1

In order to properly handle nonlocal terms in ROs it is more appropriate
to view the RO as a B\"acklund auto-transformation for
$\ell_{\mathcal{F}}(\boldsymbol{G})=0$, see \cite{gut, marv, ms03, ms-ip, p, ph} and references therein for details, and below we adhere to this point of view.
%
Likewise, an {\em adjoint recursion operator} is from this perspective a B\"acklund
auto-transformation for  $\ell_{\mathcal{F}}^{\dagger}(\boldsymbol{\gamma})=0$, cf.\ e.g.\ \cite{kv}.


As a closing remark, recall that a partial differential system is
called {\em dispersionless} (cf.\ e.g.\ \cite{ms3} and references
therein) if it can be written in the form of a quasilinear
homogeneous first-order system
\begin{equation}
\sum_{j=1}^d \sum_{\alpha=1}^n  A_{I\alpha}^j(\bi{u})\frac{\p u^\alpha}{\p x^j}=0,\quad\mbox{where}\quad I=1,\dots,m,\quad m\geq N.
\label{hds}
\end{equation}


The class of systems (\ref{hds}) is quite rich: for instance, quasilinear scalar second-order PDEs which do not explicitly involve the dependent variable $u$,
\begin{equation}\label{2or}
\sum_{i=1}^d \sum_{j=i}^d f^{ij}(\vec x,\partial u/\partial x^1,
\dots, \partial u/\partial x^d)\partial^2 u/\partial x^i\partial
x^j=0,
\end{equation}
can be brought into the form (\ref{hds}) by putting
$\bi{u}=(\partial u/\partial x^1, \dots, \partial u/\partial x^d)^T$ and hence any equation (\ref{2or}) is dispersionless.

In particular, equations (\ref{pe}), (\ref{ase}), and (\ref{equ}) considered in the examples below are of the form (\ref{2or}) and therefore dispersionless.

\section{Recursion operators}\label{rk}
Consider the following differential operators in total derivatives:
\begin{equation}\label{ablm}
\begin{array}{l}
A_i=A_i^0+\sum\limits_{j=1}^d
A_{i}^{j} D_{x^{j}},\quad B_i=B_i^0+\sum\limits_{j=1}^d B_{i}^{j} D_{x^{j}},\quad i=1,2,\\
L=L^0+\sum\limits_{k=1}^d L^{k} D_{x^{k}},\quad
M=M^0+\sum\limits_{k=1}^d M^{k} D_{x^{k}}.
\end{array}
\end{equation}
Here $A_i^j=A_i^j(\vec x, [\boldsymbol{u}])$ and $B_i^j=B_i^j(\vec x, [\boldsymbol{u}])$ for $i=1,2$ and $j=1,\dots,d$ are scalar functions,
$A_i^0=A_i^0(\vec x, [\boldsymbol{u}])$ and $B_i^0=B_i^0(\vec x, [\boldsymbol{u}])$ for $i=1,2$
are $N\times N$ matrices,
$L^{k}=L^{k} (\vec x, [\boldsymbol{u}])$ for $k=0,\dots,d$ are $N\times m$ matrices,
and $M^{k}=M^{k} (\vec x, [\boldsymbol{u}])$ for $k=0,\dots,d$
are $m\times N$ matrices.

\begin{prop}\label{Th1}Suppose that for a system (\ref{sys}) there exist the operators $A_i,B_i,L,M$
of the form (\ref{ablm}) such that
\begin{eqnarray}
\hspace*{-90mm}\mbox{i)}&&
[A_1,A_2]=0, \label{LABa}\\
\hspace*{-90mm}\mbox{ii)}&&
[B_1,B_2]=0,\label{LABb}\\
\hspace*{-90mm}\mbox{iii)}&&
(A_1 B_2-A_2 B_1)= L\circ \ell_{\mathcal{F}},
\label{LAB}\\
\hspace*{-90mm}\mbox{iv)}&&
\ell_{\mathcal{F}}=
M \circ (B_1 A_2-B_2 A_1),
\label{LAB2}\\
\hspace*{-90mm}\mbox{v)}&&
\mbox{there exist $p,q\in\{1,\dots,d\}$, $p\neq q$,
such that we can}\nonumber\\
\hspace*{-80mm}
&&
\mbox{express $D_{x^p}\tilde{\boldsymbol{U}}$ and $D_{x^q}\tilde{\boldsymbol{U}}$ from the relations}\nonumber\\
\hspace*{-80mm}&&\quad A_i (\tilde{\boldsymbol{U}})=B_i (\boldsymbol{U}),\quad i=1,2.\label{R}
\end{eqnarray}
Then relations (\ref{R}) define a recursion operator for (\ref{sys}), i.e.,
whenever $\boldsymbol{U}$ is a (possibly nonlocal)
symmetry
for (\ref{sys}), so is
$\tilde{\boldsymbol{U}}$ defined by
(\ref{R}).
\end{prop}
{\em Proof.} First of all, if $\boldsymbol{U}$ is a (possibly nonlocal)
symmetry for (\ref{sys}), then the system (\ref{R}) for $\tilde{\boldsymbol{U}}$ is compatible
by virtue of (\ref{LABa}) and (\ref{LAB}).

We now have $\ell_{\mathcal{F}}(\tilde{\boldsymbol{U}})=0$ by virtue of (\ref{LABb}) and (\ref{LAB2}),
and hence $\tilde{\boldsymbol{U}}$ is a shadow of symmetry for our system. $\Box$\looseness=-1

As an aside
note that the condition (\ref{LAB}) in a slightly
different form
has appeared 
in \cite{dun-09, fkr, gs}, and was given there a certain geometric
interpretation.\looseness=-1

In complete analogy with the above we can also readily prove the counterpart of
Proposition~\ref{Th1} for adjoint recursion operators.
\begin{prop}\label{Th2}Suppose that for a system (\ref{sys}) there exist the operators $A_i,B_i,L,M$
of the form (\ref{ablm}) but with 
$A_i^0$ and $B_i^0$ being $m \times m$ rather than $N\times N$ matrices.
Further assume that these operators are such that 
\begin{eqnarray}
\hspace*{-90mm}\mbox{i)}&& [A_1,A_2]=0, \label{comARO}\\
\hspace*{-90mm}\mbox{ii)}&& [B_1,B_2]=0,\label{comAROb}\\
\hspace*{-90mm}\mbox{iii)}&& \ell_{\mathcal{F}}^{\dagger}=
L\circ (B_1 A_2-B_2 A_1),
\label{LABARO3b}\\
\hspace*{-90mm}\mbox{iv)}&& (A_1 B_2-A_2 B_1)=M\circ \ell_{\mathcal{F}}^{\dagger},
\label{LABN2a}\\
\hspace*{-90mm}\mbox{v)}&& \mbox{there exist $p,q\in\{1,\dots,d\}$, $p\neq q$,
such that we can}\nonumber\\
\hspace*{-80mm}&& \mbox{express $D_{x^p}\tilde{\boldsymbol{\gamma}}$ and $D_{x^q}\tilde{\boldsymbol{\gamma}}$ from the relations}\nonumber\\
\hspace*{-80mm}&& A_i (\tilde{\boldsymbol{\gamma}})=B_i (\boldsymbol{\gamma}),\quad i=1,2.\label{ARO}
\end{eqnarray}

Then (\ref{ARO}) defines an adjoint recursion operator for (\ref{sys}), i.e.,
whenever $\boldsymbol{\gamma}$
is a cosymmetry for (\ref{sys}), then so is $\tilde{\boldsymbol{\gamma}}$
defined by (\ref{ARO}).\looseness=-1
\end{prop}


\section{Recursion operators from Lax pairs}\label{roltr}
We start with the following
\begin{conj}\label{c}
Under the assumptions of Proposition~\ref{Th1}
or Proposition~\ref{Th2} the operators
$\mathcal{L}_i=\lambda A_i- B_i$, $i=1,2$, where $\lambda$ is a
spectral parameter, satisfy $[\mathcal{L}_1,\mathcal{L}_2]=0$, i.e.,
they constitute a Lax pair for (\ref{sys}).
\end{conj}

While so far we were unable to prove this result in full generality,
to the best of our knowledge it holds for all examples of (adjoint) ROs
for multidimensional dispersionless systems that are known to date.

Thus, if our system (\ref{sys}) admits an (adjoint) RO, it also usually admits a
Lax pair which is linear in the spectral parameter.

Hence, a natural source
of $A_i,B_i$ and $L, M$ satisfying the conditions
of Proposition~\ref{Th1} 
is provided by the
Lax pairs for (\ref{sys}) of the form\looseness=-1
\begin{equation}\label{Lax0}
\mathcal{L}_i\psi=0,\quad i=1,2,
\end{equation}
with $\mathcal{L}_i$ linear in $\lambda$ such that $\psi$ is a
nonlocal symmetry of (\ref{sys}), i.e., we have
$\ell_\mathcal{F}(\psi)=0$, cf.\ e.g.\ \cite{p}. Then putting
$\mathcal{L}_i=\lambda B_i- A_i$ or $\mathcal{L}_i=\lambda A_i- B_i$
gives us natural candidates for $A_i$ and $B_i$ which then should be
checked against the conditions of Proposition~\ref{Th1}, and, if the
latter hold, yield a RO for (\ref{sys}).\looseness=-1

Even if it happens to be impossible to construct the recursion operator from (\ref{Lax0}),
one still can construct infinite series of nonlocal symmetries for (\ref{sys}) by expanding $\psi$ into
formal Taylor or Laurent series in $\lambda$, cf.\ e.g.\ \cite{p, as, mos}.\looseness=-1




It is important to stress that the Lax pair with the
operators $\mathcal{L}_{i}$ employed in the above construction does
not have to be the {\em original} Lax pair of our system (\ref{sys}).
In general, we should {\em custom tailor} the operators
$\mathcal{L}_{1,2}$ constituting the Lax pair for the construction
in question, so that the solutions of the associated linear problem
(\ref{Lax0}) are (in general nonlocal) symmetries, i.e., satisfy the linearized version of
our system.

The natural building blocks for these $\mathcal{L}_i$ are the
original Lax operators $\mathscr{X}_i$, their formal adjoints
$\mathscr{X}_i^\dagger$, and the operators
$\mathrm{ad}_{\mathscr{X}_i}=[\mathscr{X}_i,\cdot]$, but in general
one has to twist them, cf.\ e.g.\ Example 4 below.\looseness=-1

More explicitly, suppose that the system under study, i.e.,
(\ref{sys}), admits a Lax pair of the form
\begin{equation}\label{lax-x}
\mathscr{X}_i\psi=0,\quad i=1,2,
\end{equation}
where $\mathscr{X}_i$ are linear
in the spectral parameter $\lambda$
but $\psi$ does {\em not} satisfy $\ell_{\mathcal{F}}(\psi)=0$,
i.e., $\psi$ is not a (nonlocal) symmetry for (\ref{sys}).

Then we can seek for a nonlocal symmetry of (\ref{sys}) of the
form 
\begin{equation}\label{phi}
\Phi=\Phi(\vec x, [\boldsymbol{u},\psi,\chi,\vec{\zeta}]),
\end{equation}
i.e., for a vector function of $\vec x$, and of $\boldsymbol{u}$,
$\psi$, $\chi$, and $\vec{\zeta}=(\zeta^1,\dots,\zeta^d)^T$ and a finite number of the
derivatives of these variables.

Here $\chi$ satisfies the system \[ \mathscr{X}^\dagger_i\chi=0,\ i=1,2,\]
and $\mathscr{Z}=\sum\limits_{j=1}^d \zeta^j\partial/\partial x^j$
satisfies
\[
[\mathscr{X}_i,\mathscr{Z}]=0,\ i=1,2.
\]

Moreover, we should also require that there exist the operators
$\mathscr{L}_i$ which are linear in $\lambda$ and such that
\[\mathscr{L}_i\Phi=0,\quad i=1,2.\]
Then one should extract $A_i$ and $B_i$ for
Proposition~\ref{Th1} from these $\mathscr{L}_i$ rather than
from the original $\mathscr{X}_i$.

The study of specific examples strongly suggests that 
the nonlocal variables $\psi$, $\chi$ and $\vec{\zeta}$ usually do not mix with each other and enter $\Phi$ linearly,
so it typically suffices to look for nonlocal symmetries $\Phi$ of (\ref{sys}) in
one of the following forms,\looseness=-1
\begin{eqnarray}
\Phi^\alpha&=&\sum\limits_{s=1}^r a^\alpha_s \psi^s+
\sum\limits_{s=1}^r\sum\limits_{k=1}^d
a^{k}_s D_{x^{k}} (\psi^s),
\label{psi1}\\
\Phi^\alpha&=&\sum\limits_{s=1}^r b^{\alpha, s} \chi_s+
\sum\limits_{s=1}^r\sum\limits_{k=1}^d
b^{\alpha,s,k} D_{x^{k}}(\chi_s),
\label{psi2}\\
\Phi^\alpha&=&\sum\limits_{j=1}^d c^\alpha_j \zeta^j+
\sum\limits_{k,l=1}^d
c^{\alpha,k}_l D_{x^{k}} (\zeta^l),
\label{psi3}
\end{eqnarray}
where $\alpha=1,\dots,N$, $p\in\mathbb{N}$,
$a^\alpha_s$, $a^{\alpha,k}_s$, $b^{\alpha, s}$, $b^{\alpha,s,k}$,
$c^\alpha_j$ and $c^{\alpha,k}_l$ are local functions,
and $r$ is the number of components of $\psi$,
instead of general form (\ref{phi}).

Furthermore,
it is often possible to restrict
oneself to considering only zero-order terms in the nonlocal variables, i.e., use even simpler Ans\"atze
$\Phi^\alpha=\sum\limits_{s=1}^r a^\alpha_s \psi^s$,
$\Phi^\alpha=\sum\limits_{s=1}^r b^{\alpha, s} \chi_s$, and
$\Phi^\alpha=\sum\limits_{j=1}^d c^\alpha_j \zeta^j$ instead
of (\ref{psi1}), (\ref{psi2}) and (\ref{psi3}).\looseness=-1

A similar approach, with symmetries replaced by cosymmetries and
Proposition~\ref{Th1} by Proposition~\ref{Th2} can, of course, be
applied to the construction of adjoint recursion operators.

Finally, let us point out that the assumption of
linear dependence of the Lax operators on $\lambda$ is not as restrictive as it seems:
a great many of known today multidimensional dispersionless
hierarchies include systems with this property and conversely,
if a system admits an RO of the form described
in Proposition~\ref{Th1} then it usually possesses a Lax-type
representation constructed as per Conjecture~\ref{c} whose operators are linear in $\lambda$.\looseness=-1

Moreover, in some cases the nonlinearity in $\lambda$ can be removed 
upon a proper rewriting of the Lax pair.

Consider for example the Pavlov equation \cite{mvp}
\begin{equation}\label{pe}
u_{yy}+u_{xt}+u_x u_{xy}-u_y u_{xx}=0,
\end{equation}
so in our notation $d=3$, $N=1$, $x^1=x$, $x^2=y$, $x^3=t$, $u^1=u$.

Equation (\ref{pe}) possesses a Lax pair of the form \cite{mvp} (cf.\ also \cite{mik, m-a-s, ms})
\begin{equation}\label{lax-pe}
\psi_y=(-u_x-\lambda) \psi_x,\quad \psi_t=(-\lambda^2-\lambda u_x+u_y) \psi_x.
\end{equation}

This Lax pair is quadratic in $\lambda$ but rewriting its second equation as
\[
\psi_t=\lambda(-\lambda- u_x)\psi_x+u_y\psi_x
\]
and substituting $\psi_y$ for $(-\lambda- u_x)\psi_x$, we obtain
\[
\psi_t=\lambda\psi_y+u_y\psi_x,
\]
and thus (\ref{lax-pe}) can be rewritten in the form linear in $\lambda$:
\begin{equation}\label{lax-pe1}
\psi_y=(-u_x-\lambda) \psi_x,\quad \psi_t=-\lambda\psi_y+u_y\psi_x.
\end{equation}

\section{Examples}\label{ex}
\subsection*{Example 1} We now continue the study of the Lax pair (\ref{lax-pe1}) for the Pavlov equation. It is readily checked that if $\psi$ satisfies (\ref{lax-pe1})
then $\Phi=1/\psi_x$ is a nonlocal symmetry for
(\ref{pe}), i.e., $U=1/\psi_x$ satisfies the
linearized version of (\ref{pe}), 
\begin{equation}\label{lpe}
U_{yy}+U_{xt}+u_x U_{xy}-u_y U_{xx}+ u_{xy} U_x- u_{xx}U_y=0.
\end{equation}

Now by virtue of (\ref{lax-pe1}) the function $\Phi=1/\psi_x$ also satisfies a pair of linear equations of the form $\mathscr{L}_i\Phi=0$, $i=1,2$, where
\[
\mathscr{L}_1=-D_y+(\lambda-u_x)D_x+u_{xx},\quad \mathscr{L}_2=D_t+\lambda D_y-u_y D_x+u_{xy}.
\]
Now, as $\Phi=1/\psi_x$ satisfies (\ref{lpe})
we can identify $\mathscr{L}_i$ with $\mathcal{L}_i$
from (\ref{Lax0}) and put $\mathscr{L}_i=\lambda A_i- B_i$, where
\[
A_1=D_x, \quad A_2=D_y, \quad B_1=D_y+u_x D_x-u_{xx},\quad
B_2=-D_t+u_y D_x-u_{xy}.\]
Then it is easily seen that these operators satisfy
the conditions of Proposition~\ref{Th1} for suitably chosen $L$ and
$M$, so we arrive at the recursion operator for the Pavlov equation
given by the formulas
\[
\tilde U_x=U_y+u_x U_x-u_{xx} U,\quad \tilde U_y=-U_t+u_y U_x-u_{xy} U,
\]
which maps a (possibly nonlocal) symmetry $U$ to a new nonlocal symmetry
$\tilde U$. This is nothing but the recursion operator found in
\cite{ms} rewritten as a B\"acklund auto-transformation for (\ref{lpe}).

\subsection*{Example 2} Consider the equation \cite{sch, sch2}
\begin{equation}\label{ghe}
a u_{xy}u_{zt} + b u_{xz}u_{yt}+c u_{xt}u_{yz} = 0,\qquad a+b+c=0,
\end{equation}
where $a,b,c$ are constants.
We have $d=4$, $x^1=x$, $x^2=y$, $x^3=z$, $x^4=t$, $N=1$, $u^1=u$.

In \cite{df}, where equation (\ref{ghe}) was rediscovered in the course of classification of integrable symplectic Monge--Amp\`ere equations, thos equation was named the general heavenly equation, and it was shown how to rewrite (\ref{ghe}) in the standard dispersionless form (\ref{hds}). Note that (\ref{ghe}) describes {\em
inter alia} \cite{sch, sch2, mas} a class of anti-self-dual solutions of the Einstein
field equations.

By definition, a (nonlocal) symmetry $U$ of (\ref{ghe}) satisfies
the linearized version of (\ref{ghe}), that is,
\begin{equation}\label{linghe}
a u_{xy}U_{zt} +a u_{zt}U_{xy} + b u_{xz}U_{yt}+
b u_{yt}U_{xz}+c u_{xt}U_{yz} +c u_{yz}U_{xt}= 0. 
\end{equation}

Equation (\ref{ghe}) admits \cite{sch,sch2,df,mas}
a Lax pair with the operators
\begin{equation}\label{lax}
\begin{array}{l}
L_1=(1+c\lambda)D_x-\displaystyle\frac{u_{xz}}{u_{zt}}D_t-c\lambda \displaystyle\frac{u_{xt}}{u_{zt}}D_z,\\[5mm]
L_2=(1-b\lambda)D_y-\displaystyle\frac{u_{yz}}{u_{zt}}D_t+b\lambda \displaystyle\frac{u_{yt}}{u_{zt}}D_z.
\end{array}
\end{equation}
It is readily seen that any $\psi$ that satisfies $L_i\psi=0$, $i=1,2$, also satisfies (\ref{linghe}), i.e., it is a nonlocal symmetry for (\ref{ghe}), so in spirit of Conjecture~\ref{c} 
let
\[
\begin{array}{l}
A_1=c D_x-c \displaystyle\frac{u_{xt}}{u_{zt}}D_z,\quad
A_2=-b D_y+b \displaystyle\frac{u_{yt}}{u_{zt}}D_z,\\[5mm]
B_1= -D_x+\displaystyle\frac{u_{xz}}{u_{zt}}D_t,\quad
B_2= -D_y+\displaystyle\frac{u_{yz}}{u_{zt}}D_t.
\end{array}
\]
Then all conditions of Proposition~\ref{Th1} are readily seen to be
satisfied, e.g.\ we have $B_1 A_2-B_2 A_1=(1/u_{zt})\ell_F$, where
$F=a u_{xy}u_{zt} + b u_{xz}u_{yt}+c u_{xt}u_{yz}$.

Hence 
the relations
\begin{equation}\label{ghero}
\displaystyle\tilde{U}_x=\frac{u_{xz}U_t+c u_{xt}\tilde{U}_z
-u_{zt} U_x}{cu_{zt}},\
\displaystyle\tilde{U}_y=-\frac{u_{yz}U_t-b u_{yt}\tilde{U}_z-u_{zt} U_y}{bu_{zt}},
\end{equation}
where $U$ is a (possibly nonlocal) symmetry of (\ref{ghe}), define a new symmetry $\tilde U$
and thus a novel recursion operator (\ref{ghero}) for (\ref{ghe}), i.e., a B\"acklund auto-transformation for (\ref{linghe}). 
This operator has first appeared in the preprint version of the present paper (arXiv:1501.01955v1) and was later rediscovered, in a slightly different form, in \cite{msy}.


In closing note that there exists a transformation \cite{sch, sch2} relating (\ref{ghe}) and the first heavenly equation, so in principle it should be possible to employ the transformation in question to obtain the recursion operator (\ref{ghero}) from a suitable combination of ROs for the first heavenly equation found in \cite{sy}. However, this transformation is very complicated, and applying our approach provides a much easier way to find the recursion operator in question.\looseness=-1

\subsection*{Example 3} Consider the Mart\'\i{}nez Alonso--Shabat \cite{m-a-s} equation
\begin{equation}\label{ase}
u_{yt}=u_z u_{xy}-u_y u_{xz}
\end{equation}
so now $d=4$, $N=1$, $x^1=x$, $x^2=y$, $x^3=z$, $x^4=t$, $u^1=u$.

Eq.(\ref{ghe}) admits \cite{m3} a Lax pair with the
operators
\begin{equation}\label{lax-ase}
L_1=D_y-\lambda u_y D_x,\quad 
L_2=D_z-\lambda u_z D_x+\lambda D_t.
\end{equation}
Consider the action of operators
$\mathcal{L}_i=\mathrm{ad}_{L_i}=[L_i,\cdot]$
on operators of the form $w D_x$ and spell it out as
\[
\mathrm{ad}_{L_i}(w D_x)\equiv((\lambda B_i-A_i)w)D_x,
\]
whence $A_1= D_y$, $A_2= D_z$, $B_1=u_y D_x-u_{xy}$,
$B_2=u_z D_x-D_t-u_{xz}$.

All conditions of Proposition~\ref{Th1} are again satisfied, so the relations
\begin{equation}\label{asero}
\displaystyle\tilde{U}_y=u_y U_x-u_{xy}U,\quad
\displaystyle\tilde{U}_z=u_z U_x-u_{xz}U-U_t,
\end{equation}
where $U$ is any (possibly nonlocal) symmetry of (\ref{ase}), define
a new symmetry $\tilde U$, i.e., Eq.(\ref{asero}) provides a recursion
operator for (\ref{ase}). This operator was already found in
\cite{m3}.

Applying (\ref{asero}) to the symmetry $U=u_x$ we obtain (modulo an arbitrary function of $x$ and $t$ resulting from the integration) a nonlocal symmetry of the form  $\tilde{U}=w-u_x^2/2$, where the nonlocal variable $w$ is defined by the relations\looseness=-1
\[
\begin{array}{l}
\displaystyle w_y=u_y u_{xx},\qquad
\displaystyle w_z=u_z u_{xx}-u_{xt}.
\end{array}
\]
This example suggests that the recursion operator (\ref{asero}) can be further simplified: a new symmetry $\tilde{U}$ can be constructed by putting
\[\tilde{U}=W-u_x U,
\]
where $W$ is defined by the relations
\[
\displaystyle W_y=u_y U_x+u_{x}U_y,\quad
\displaystyle W_z=u_z U_x+u_{x}U_z-U_t.
\]
Thus, the recursion operators produced within our approach are not
necessarily in the simplest possible form.

\subsection*{Example 4} Consider a system \cite{fk, fkk}
\begin{equation}
m_{t}=n_{x}+nm_{r}-mn_{r},\qquad n_{z}=m_{y}+m n_{s}-n m_{s},
\label{fks}
\end{equation}
where $x,y,z,r,s,t$ are independent variables and $m,n$ are dependent variables, so
$m=N=2$, $d=6$, $\boldsymbol{u}=(m,n)^T$ and $\vec x=(x,y,z,r,s,t)^T$.

System (\ref{fks}) is integrable, as it admits a Lax pair. Namely, it can \cite{fk,fkk} be written as a condition of commutativity for the following pair of vector fields
involving a spectral parameter $\lambda$:
\be\label{cc0}
[D_{z}-m D_{s}-\lambda D_{x}+\lambda m D_{r},
D_{y}-n D_{s}-\lambda D_{t}+\lambda n D_{r}]=0.
\ee

By virtue of (\ref{fks}) there exists a potential $u$ such that
\be\label{sub}
m=u_z/u_s,\quad n=u_y/u_s.
\ee
Substituting this into (\ref{fks}) gives a single second-order equation for $u$,
\be\label{equ}
u_s u_{zt}-u_z u_{st}- u_s u_{xy}+u_y u_{sx}-u_y u_{rz}+u_z u_{ry}=0,
\ee
which can be written as a commutativity condition $[L_1,L_2]=0$ for
\[
L_1=\displaystyle D_{z}-\frac{u_z}{u_s}D_{s}-\lambda D_{x}+\lambda \frac{u_z}{u_s} D_{r},\
L_2=\displaystyle D_{y}-\frac{u_y}{u_s}D_{s}-\lambda D_{t}+\lambda \frac{u_y}{u_s}D_{r}.
\]

Let $\chi$ satisfy $L_i^\dag \chi=0$, $i=1,2$.
Then $\zeta=u_s/\chi$ is a nonlocal symmetry for (\ref{equ}),
and we can readily obtain a linear system for $\zeta$ of the form
$\mathcal{L}_i\zeta=0$, $i=1,2$, from that for $\chi$.\looseness=-1

Upon 
spelling out $\mathcal{L}_i$ as $\mathcal{L}_i=A_i-\lambda B_i$
and checking that the conditions of Theorem~\ref{Th1} are satisfied
for suitable $L$ and $M$ we arrive at
the following novel RO 
for (\ref{equ}):
\be\label{invro1}
\ba{l}
\tilde{U}_y =\ds \frac{u_y}{u_s} \tilde{U}_s-\frac{u_y}{u_s} U_r+ U_t-\frac{(u_{st}-u_{ry})}{u_s} U, \\[7mm]
\tilde{U}_z =\ds \frac{u_z}{u_s} \tilde{U}_s-\frac{u_z}{u_s} U_r+ U_x+\frac{(u_{rz}-u_{sx})}{u_s} U.
\ea
\ee
Here again $U$ is a (possibly nonlocal) symmetry of (\ref{equ}), and $\tilde U$
is a new (again possibly nonlocal) symmetry for (\ref{equ}).

\subsection*{Example 5} Consider (see e.g.\ \cite{ac, atcmp, Mason} and references therein)
the 
Lax operators
\begin{equation}\label{kmlax}
L_1=D_y+\partial_{\tilde x} K-\lambda D_{\tilde x},\quad
L_2=D_x+\partial_{\tilde y} K-\lambda D_{\tilde y}.
\end{equation}
Here $K$ takes values in a (matrix) Lie algebra $\mathfrak{g}$.\looseness=-2

The commutativity condition
$[L_1,L_2]=0$
yields the equation for $K$
\begin{equation}\label{kme}
\partial_x \partial_{\tilde x} K-\partial_y \partial_{\tilde y} K
-[\partial_{\tilde x} K,\partial_{\tilde y} K]=0.
\end{equation}

Note that (\ref{kme}) can be rewritten in the dispersionless form of the type (\ref{hds})
upon introducing (cf.\ e.g.\ \cite{Mason, ph}) an additional dependent variable $J$
taking values in the (matrix) Lie group $G=\exp(\mathfrak{g})$.
Namely, the dispersionless form of (\ref{kme}) reads
\begin{equation}\label{jk}
\partial_x J= J \partial_{\tilde y} K,\quad \partial_y J= J \partial_{\tilde x} K.
\end{equation}
The compatibility condition $\partial_x\partial_y J=\partial_y\partial_x J$ yields an equation
\begin{equation}\label{jme}
\partial_{\tilde x} (J^{-1}\partial_x J) = \partial_{\tilde y} (J^{-1}\partial_y J),
\end{equation}
while (\ref{kme}) is identically satisfied by virtue of (\ref{jk}).
In other words, (\ref{jk}) defines a B\"acklund transformation relating (\ref{kme}) and (\ref{jme}).

The gauge potentials defined by
setting $A_x=\partial_{\tilde y} K$, $A_y=\partial_{\tilde x} K$, $A_{\tilde x} =0$, $A_{\tilde y}=0$,
satisfy the anti-self-dual Yang--Mills equations on $\mathbb{R}^4$ with Euclidean metric $g=\mathrm{diag}(+1,+1,+1,+1)$ in the coordinates $X^i$
such that $\sqrt{2}x=X^1+ \mathrm{i} X^2$, $\sqrt{2}\tilde x=X^1- \mathrm{i} X^2$, $\sqrt{2}y=X^3- \mathrm{i} X^4$,
$\sqrt{2}\tilde y=X^3+ \mathrm{i} X^4$, where $\mathrm{i}=\sqrt{-1}$, cf.\ e.g.
\cite{Mason}. In this context $K$ is known as the Yang $K$-matrix. Conversely, any sufficiently smooth solution of
the anti-self-dual Yang--Mills equations can, up to a suitable gauge transformation, be obtained in this way \cite{Mason}.\looseness=-1

Thus, the anti-self-dual Yang--Mills equations, which play an important role in modern physics, cf.\ e.g.\ \cite{Mason} and references therein, can, in a suitable gauge, be written in dispersionless form, exactly as stated in Introduction.\looseness=-1

It is natural to assume that in the Lax pair equations $L_i\psi=0$ we
have $\psi$ taking values in $G$, but then $\psi$ cannot be a symmetry
for (\ref{kme}), since a symmetry must live in $\mathfrak{g}$ just like $K$.
However, we can get a Lax pair
for $\Phi\in\mathfrak{g}$ using the adjoint action: $[L_1,\Phi]=[L_2, \Phi]=0$, i.e., we use the Lax operators $\mathcal{L}_i=\mathrm{ad}\ L_i$.

Then all conditions of Proposition~\ref{Th1} are satisfied for $A_i$ and $B_i$ obtained
by spelling out $\mathcal{L}_i$ as
$\mathcal{L}_i=\lambda A_i -B_i$.
This yields an RO for (\ref{kme}) of the form
\begin{equation}\label{kmero}
V_{\tilde x}=U_y+[\partial_{\tilde x} K,U],\quad V_{\tilde y}=U_x+[\partial_{\tilde y} K,U],
\end{equation}
which produces a symmetry $V$ for (\ref{kme}) from another symmetry $U$.
This RO is known from the literature, see \cite{blr, ph, atcmp}.

If in the above construction we take $\mathfrak{g}=\mathfrak{diff}(\mathbb{R}^N)$, i.e., we put
$K=\sum\limits_{i=1}^N u_i
\partial/\partial z^i,$
where
$u_i=u_i(x,y,\tilde x,\tilde y,\allowbreak z^1,\dots,\allowbreak z^N)$ are scalar functions,
then (\ref{kme}) becomes \cite{ms3} the Manakov--Santini system \cite{ms-pla},
an integrable system in $(N+4)$ independent variables, and
(\ref{kmero}) provides a recursion operator for this system which was
found earlier in \cite{ms-ip} (see also \cite{ms3}) by a different
method.

Further examples of recursion operators obtained using our approach
can be found e.g.\ in \cite{km}. 

\section{Concluding remarks}\label{co}

We presented above a method of construction for recursion operators
and adjoint recursion operators for a broad class of
multidimensional integrable systems that can be written as
commutativity conditions for a pair of vector fields, or, even more
broadly, of linear combinations of vector fields with zero-order
matrix operators, linear in the spectral parameter and free of
derivatives in the latter. This method is based on constructing a
special Lax pair for the system under study from the
original Lax pair, see Section~\ref{roltr} for
details. Note that if a system under study admits several
essentially distinct Lax pairs then our approach applied
to these Lax pairs in principle could give rise to several
essentially distinct recursion operators.

To the best of author's knowledge, the method in question
works for all known examples of multidimensional dispersionless integrable
systems admitting recursion operators in the form of B\"acklund
auto-transformation for linearized systems, including e.g.\ the $ABC$
equation, see \cite{z, ksm} and references therein,
\[
a u_x u_{yt} + b u_y u_{xt} + c u_t u_{xy} = 0,\quad a+b+c=0,
\]
the simplest (2+1)-dimensional equation of the so-called universal hierarchy
of Mart\'\i{}nez Alonso and Shabat \cite{m-a-s},
\[
u_{yy}-u_y u_{tx}+u_x u_{ty}=0,
\]
the complex Monge--Amp\`ere equation, see \cite{nsky} for its
recursion operator, first, second and modified heavenly equations,
etc. In particular, this method enabled us to find hitherto unknown
recursion operators for the general heavenly equation (\ref{ghe})
and for (\ref{equ}).
Note that our approach, when applicable, is in general computationally less demanding
than those of \cite{ms-ip} and \cite{m2}.
\looseness=-1

While all examples to which we have applied our approach so far are dispersionless,
Propositions~\ref{Th1} and \ref{Th2} and the method described in Section~\ref{roltr}
do not explicitly require this to be the case, so
it would be interesting to find out whether it is possible to extend
our approach, perhaps also letting the operators $A_i,B_i,L$ and $M$ be of degree higher than one in the process, to the construction of recursion operators
to other classes of integrable systems.

\subsection*{Acknowledgments}
The author is pleased to thank I.S. Krasil'shchik, O.I. Morozov, M.V. Pavlov and R.O. Popovych for stimulating discussions, to W.K. Schief for pointing out the references \cite{sch,sch2}, and to the anonymous referee for valuable recommendations.\looseness=-1

This research was supported in part by the Ministry of Education, Youth and Sports of the Czech Republic (M\v{S}MT \v{C}R) under RVO funding for I\v{C}47813059, and by the Grant Agency of the Czech Republic (GA \v{C}R) under grant P201/12/G028.



\begin{thebibliography}{99}
\footnotesize
\itemsep=-0.05mm
\renewcommand{\baselinestretch}{1}

\bibitem{ac}M.J. Ablowitz, P.A. Clarkson, Solitons, nonlinear evolution equations and inverse scattering, Cambridge University
Press, Cambridge, 1991.

\bibitem{atcmp} M. Ablowitz, S. Chakravarty, L.A. Takhtajan, A self-dual Yang-Mills hierarchy and its reductions to integrable systems in $1+1$ and $2+1$ dimensions, Comm. Math.
Phys. 158 (1993), no. 2, 289--314.

\bibitem{a-s}A.A. Agrachev, Yu.L. Sachkov,
Control theory from the geometric viewpoint,
Berlin, Springer, 2004. 

\bibitem{bp}A. Bihlo, R.O. Popovych, Group classification of linear evolution equations,
J. Math. Anal. Appl. 448 (2017), no. 2, 982--1005, arXiv:1605.09251


\bibitem{bl}M. B\l aszak, Multi-Hamiltonian Dynamical Systems, Berlin, Springer, 1998.


\bibitem{blu} G.W. Bluman, A.F. Cheviakov, S.C. Anco, Applications of symmetry methods to partial differential equations, 
Springer, New York, 2010.

\bibitem{blr} M. Bruschi, D. Levi and O. Ragnisco, Nonlinear partial differential equations and B\"acklund transformations related to the 4-dimensional self-dual Yang-Mills equations, Lett. Nuovo Cimento 33 (1982), 263--266.

\bibitem{ca}F. Calogero, Why are certain nonlinear PDEs both widely applicable and
integrable?, in {\em What is integrability?}, Springer, Berlin, 1991, 1--62.


\bibitem{cf}D. Catalano Ferraioli, L.A. de Oliveira Silva,  Local isometric immersions of pseudo\-spherical surfaces described by evolution equations in conservation law form, J. Math. Anal. Appl. 446  (2017),  no. 2, 1606--1631.\looseness=-1

\bibitem{cgs}J. Cie\'sli\'nski, P. Goldstein, A. Sym, Isothermic surfaces in $\mathbf{E}^3$ as soliton surfaces, Phys. Lett. A  205  (1995),  no. 1, 37--43, arXiv:solv-int/9502004


\bibitem{df}B. Doubrov and E.V. Ferapontov, On the integrability of symplectic
Monge--Amp\`ere equations, J. Geom. Phys. 60 (2010), 1604--1616,
arXiv:0910.3407


\bibitem{dun-09}M. Dunajski,
The twisted photon associated to hyper-Hermitian four-manifolds, J.
Geom. Phys. 30 (1999), no. 3, 266--281.

\bibitem{dun}M. Dunajski, Solitons, Instantons and Twistors, Oxford University
Press, 2010.


\bibitem{fk}E.V. Ferapontov,
K.R. Khusnutdinova, Hydrodynamic reductions of multidimensional
dispersionless PDEs: the test for integrability, J. Math. Phys. 45
(2004), 2365--2377, arXiv:nlin/0312015

\bibitem{fkk}E.V. Ferapontov,
K.R. Khusnutdinova, C. Klein, On linear degeneracy of integrable quasilinear systems in higher dimensions,
Lett. Math. Phys. 96 (2011), no. 1-3, 5--35, arXiv:0909.5685.

\bibitem{fkr} E.V. Ferapontov, B. Kruglikov, Dispersionless integrable systems in 3D and Ein\-ste\-in--Weyl geometry,
J. Differential Geom. 97 (2014), no. 2, 215--254, arXiv:1208.2748 \looseness=-1 %

\bibitem{gs}J.D.E. Grant, I.A.B. Strachan, Hypercomplex integrable systems, Nonlinearity 12 (1999) 1247--1261, arXiv:solv-int/9808019

\bibitem{gut}G.A. Guthrie, Recursion operators and non-local symmetries, Proc. R. Soc. 
A 446 (1994), no. 1926, 107--114.

\bibitem{kkv}P. Kersten, I.S. Krasil'shchik, A. Verbovetsky, R.  Vitolo,  Hamiltonian structures for general PDEs, in {\em Differential equations: geometry, symmetries and integrability},  187--198,
Springer, Berlin, 2009, arXiv:0812.4895.

\bibitem{ksm}I.S. Krasil'shchik, A. Sergyeyev, O.I. Morozov, Infinitely many nonlocal conservation laws for the $ABC$ equation with $A+B+C\neq 0$, Calc. Var. Partial Differential Equations 55 (2016), no. 5, article 123,  arXiv:1511.09430

\bibitem{kv}J. Krasil'shchik, A. Verbovetsky, Geometry of jet spaces and integrable systems,
J. Geom. Phys. 61 (2011), no. 9, 1633--1674, arXiv:1002.0077

\bibitem{km}B. Kruglikov, O. Morozov, Integrable dispersionless PDEs in 4D, their symmetry pseudogroups and deformations, Lett. Math. Phys. 105 (2015), no. 12, 1703--1723, arXiv:1410.7104.

\bibitem{l}J. Lenells, Exactly Solvable Model for Nonlinear Pulse Propagation in Optical Fibers,  Stud. Appl. Math. 123 (2009), no. 2, 215--232.


\bibitem{ma}W.X. Ma, R.K. Bullough, P.J. Caudrey, W.I. Fushchych, Time-dependent symmetries of variable-coefficient evolution equations and graded Lie algebras, J. Phys. A: Math. Gen. 30 (1997), no. 14, 5141--5149.

\bibitem{mns} A.A. Malykh, Y. Nutku, M.B. Sheftel, Partner symmetries of the
complex Monge-Amp\`ere equation yield hyper-K\"ahler metrics without
continuous symmetries, J. Phys. A: Math. Gen. 36 (2003), 10023--10037.

\bibitem{m-s}A.A. Malykh, Y. Nutku, M.B. Sheftel, Partner symmetries and non-invariant solutions of four-dimensional heavenly equations,
J. Phys. A: Math. Gen. 37 (2004), 7527--7545, arXiv:math-ph/0403020

\bibitem{mas}A.A. Malykh, M.B. Sheftel, General heavenly equation governs anti-self-dual gravity,
J. Phys. A: Math. Theor. 44 (2011), no. 15, 155201, 11 pp, arXiv:1011.2479

\bibitem{ms-pla}S.V. Manakov and P.M. Santini,
Inverse scattering problem for vector fields and the Cauchy problem
for the heavenly equation, {\it Phys. Lett. A} 359 (2006),
no. 6, 613--619, arXiv:nlin/0604024.

\bibitem{ms}S.V. Manakov, P.M. Santini,
A hierarchy of integrable partial differential equations in 2+1 dimensions associated with one-parameter families of one-dimensional vector fields, Theor. Math. Phys. 152 (2007), no. 1, 1004--1011, arXiv:nlin/0611047.

\bibitem{ms2}S.V. Manakov, P.M. Santini, On the solutions of the second heavenly and Pavlov equations, J. Phys. A: Math. Theor. 42 (2009), no. 40, article 40401, arXiv:0812.3323.

\bibitem{ms3}S.V. Manakov, P.M. Santini, Integrable dispersionless PDEs arising as commutation condition of pairs of vector fields, J. Phys. Conf. Ser. 482 (2014), paper 012029.

\bibitem{m-a-s}L. Mart\'\i{}nez Alonso, A.B. Shabat, Hydrodynamic reductions and solutions of a universal hierarchy,
Theor. Math. Phys. 104 (2004), 1073--1085


\bibitem{marv}M. Marvan, Another look on recursion operators, in {\em Differential geometry and applications (Brno, 1995)}, 393--402, Masaryk Univ., Brno, 1996.

\bibitem{ms03}M. Marvan, A. Sergyeyev, Recursion operator for the stationary Nizhnik--Veselov--Novikov equation, J. Phys. A: Math. Theor. 36 (2003), no. 5, L87--L92, arXiv:nlin/0210028.

\bibitem{ms-ip}M. Marvan, A. Sergyeyev, Recursion operators for dispersionless integrable systems in any dimension, Inverse Problems 28 (2012), paper 025011, 12 p., arXiv:1107.0784.\looseness=-1

\bibitem{Mason}L. Mason, The anti-self-dual Yang--Mills equations and their reductions, in {\em Geometry and Integrability}, eds. L. Mason and Y.
Nutku,  60--88,
Cambridge Univ.\ Press, Cambridge, 2003.\looseness=-1

\bibitem{my}A.V. Mikhailov, R.I. Yamilov, Towards classification of (2+1)-dimensional integrable
equations. Integrability conditions. I, J. Phys. A: Math. Gen.\  31  (1998),
6707--6715.\looseness=-1

\bibitem{mik}V.G. Mikhalev, Hamiltonian formalism of Korteweg--de Vries-type hierarchies,
Funct. Anal. Appl. 26 (1992), no. 2, 140--142.

\bibitem{m1}O.I. Morozov, Recursion Operators and Nonlocal Symmetries for Integrable rmdKP and rdDym Equations, preprint arXiv:1202.2308

\bibitem{m2}O.I. Morozov, A Recursion Operator for the Universal Hierarchy Equation via Cartan's Method of Equivalence,
Cent. Eur. J. Math. 12 (2014), 271--283, arXiv:1205.5748

\bibitem{m3}O.I. Morozov, The four-dimensional Mart\'\i{}nez Alonso--Shabat equation: differential coverings and recursion operators,
J. Geom. Phys. 85 (2014), 75--80, arXiv:1309.4993


\bibitem{mos}O.I. Morozov, A. Sergyeyev, The four-dimensional Mart\'\i{}nez Alonso--Shabat equation:
reductions and nonlocal symmetries, J. Geom. Phys. 85 (2014), 40--45, arXiv:1401.7942.

\bibitem{n}A. Newell, Solitons in Mathematics and Physics, SIAM, Philadelphia, 1985.

\bibitem{nd}J.M. Nunes da Costa, P.A. Damianou, Toda systems and exponents of simple Lie groups,  Bull. Sci. Math. 125 (2001), no. 1, 49--69.

\bibitem{nsky}Y. Nutku, M.B. Sheftel, J. Kalayci and D. Yaz\i{}c\i{}, Self-dual gravity is completely
integrable, J. Phys. A: Math. Theor. 41 (2008), paper 395206,
arXiv:0802.2203

\bibitem{olver}P.J. Olver, Evolution equations possessing infinitely many symmetries, J. Math. Phys. 18 (1977), no. 6, 1212--1215.

\bibitem{olverbook}P.J. Olver,  Applications of Lie groups to differential equations, 2nd ed.,
Springer-Verlag, New York, 1993.

\bibitem{p}C.J. Papachristou, Symmetry, Conserved Charges, and Lax Representations of Nonlinear Field Equations: A Unified Approach,
Electron. J. Theor. Phys. 7, No. 23 (2010) 1--16, arXiv:1001.0872

\bibitem{ph} C.J. Papachristou, B.K. Harrison, Local currents for the $GL(N,\mathbb{C})$ self-dual Yang-Mills equation,
Lett. Math. Phys. 17 (1989), no. 4, 285--288.

\bibitem{mvp}M.V. Pavlov, Integrable hydrodynamic chains, J. Math. Phys. 44
(2003), no. 9, 4134--4156, arXiv:nlin/0301010.

\bibitem{sch}W.K. Schief, Self-dual Einstein spaces via a permutability theorem for the Tzitzeica equation, Phys. Lett. A 223 (1996), no. 1-2, 55--62.

\bibitem{sch2}W.K. Schief, Self-dual Einstein spaces and a discrete Tzitzeica equation. A per\-mutability theorem link, in
{\em Symmetry and integrability of difference equations},
137--148, Cambridge University Press, Cambridge, 1999.

\bibitem{as}A. Sergyeyev, Infinite hierarchies of nonlocal symmetries of the Chen--Kontsevich--Schwarz type for the oriented associativity
equations, J. Phys. A: Math. Theor. 42 (2009), no. 40, article 404017, 15
pp., arXiv:0804.2020

\bibitem{s}A. Sergyeyev, A new class of (3+1)-dimensional integrable systems related to contact geometry, arXiv:1401.2122

\bibitem{sv}A. Sergyeyev, R. Vitolo, Symmetries and conservation laws for the Karczewska--Rozmej--Rutkowski--Infeld equation, Nonlin.\ Analysis: Real World Appl.\  32  (2016), 1--9, arXiv:1511.03975.\looseness=-1

\bibitem{sy} M.B. Sheftel, D. Yaz\i{}c\i{}, Recursion Operators and Tri-Hamiltonian
Structure of the First Heavenly Equation of Pleba\'nski, SIGMA 12
(2016), 091, 17 pp., arXiv:1605.07770

\bibitem{msy}M.B. Sheftel, A.A. Malykh, D. Yaz\i{}c\i{}, Recursion operators and bi-Hamiltonian structure of the general heavenly
equation, J. Geom. Phys. 116 (2017), 124--139, arXiv:1510.03666


\bibitem{vps}O.O. Vaneeva, R.O. Popovych, C. Sophocleous,  Equivalence transformations in the study of integrability, Phys. Scr. 89 (2014), 038003, 9 pp., arXiv:1308.5126.


\bibitem{z}I. Zakharevich,
Nonlinear wave equation, nonlinear Riemann problem, and the twistor
transform of Veronese webs, preprint arXiv:math-ph/0006001

\end{thebibliography}
\end{document}